\documentclass[final]{siamart171218}

\usepackage{amsmath}
\usepackage{amsfonts}
\usepackage{amssymb}
\usepackage{graphicx}
\usepackage{natbib}
\usepackage{url}
\usepackage{booktabs}
\usepackage{subcaption}
\usepackage{float}

\headers{Market Dynamical Systems}{A. Komarla and M. Hill}

\title{Market Dynamical Systems\thanks{Presented at the Institute for Operations Research and Management Science (INFORMS) Annual Meeting in 2024. Submitted to SIAM Journal on Applied Mathematics.}}

\author{Aparna Komarla\thanks{Data Scientist, Intel Foundry Markets and Intelligence Group (\email{aparna.komarla@intel.com})}.
\and Max Hill\thanks{Data Scientist, Intel Corporate Strategy Group (\email{max.hill@intel.com}).}}

\begin{document}

\maketitle

\begin{abstract}
We present a novel approach to modeling market dynamics using ordinary differential equations that explicitly incorporate product competitiveness and consumer behavior. Our framework treats market segments as interacting populations within a dynamical system, analogous to predator-prey models in ecology. Competitive advantages drive market share transitions through mechanistic modeling of three key processes: new product adoption, product refresh cycles, and obsolescence dynamics.
\end{abstract}



\section{Introduction}

Market dynamics present a fundamental challenge for business planning. Traditional approaches struggle to capture the complex interactions among competitive forces, consumer preferences, and technological change that drive market evolution. As \citet{arthur1997economy} and \citet{farmer2009economy} have demonstrated, traditional econometric approaches, despite their widespread adoption, fail to capture the underlying mechanisms that drive market behavior, limiting their utility for strategic decision-making. We propose an alternative approach to market modeling based on ordinary differential equations (ODEs) and constrained parameter optimization. Our framework draws inspiration from successful applications of dynamical systems theory in population biology, epidemiology, and fluid mechanics.

Markets are complex dynamical systems influenced by product competitiveness, supply chain constraints, consumer preferences, and regulatory changes. \citet{box1976time} and \citet{sims1980macroeconomics} show that for decades, autoregressive (AR) models---including ARIMA and VAR---have dominated market forecasting. However, AR models provide limited mechanistic insight because they represent markets as weighted linear combinations of historical trends. While AR models excel at pattern recognition and short-term forecasting, \citet{granger1969investigating} and \citet{hamilton1994time} show that they provide little understanding of causal mechanisms, limiting their utility for scenario analysis and policy intervention assessment.

In contrast, dynamical systems approaches explicitly model the mechanisms that drive market outcomes. As \citet{strogatz2014nonlinear} and \citet{medio2001nonlinear} demonstrate, these frameworks quantify the relationships among market forces and enable mechanistic analysis. Moreover, \citet{pearl2009causality} argues that high stakes strategic decisions require understanding how specific interventions will affect outcomes---precisely the type of prediction that linear AR models, lacking causal structure, cannot provide.

\section{Background}

\subsection{Limitations of Traditional Market Modeling Approaches}

Traditional econometric models based on autoregressive approaches have dominated market analysis for decades, valued for their mathematical tractability and short-term forecasting accuracy. As described by \citet{box1976time} and \citet{lutkepohl2005new}, ARIMA (Autoregressive Integrated Moving Average) models and their multivariate extensions such as Vector Autoregression (VAR) models represent markets as linear combinations of lagged variables, assuming that future values can be predicted based on weighted averages of historical observations.

However, these approaches suffer from several fundamental limitations that have been extensively documented in the literature:

\begin{enumerate}
    \item \textbf{Lack of Mechanistic Insight}: \citet{granger1969investigating} and \citet{hendry1995dynamic} show that AR models are sophisticated curve-fitting exercises that identify statistical patterns but fail to explain the underlying economic mechanisms driving market behavior. \citet{sims1980macroeconomics} demonstrates that these models cannot distinguish correlation from causation, rendering them unsuitable for policy analysis.
    \item \textbf{Limited Scenario Analysis}: Because these models rely on historical correlations, they cannot predict responses to unprecedented interventions. \citet{lucas1976econometric} showed in his famous critique that structural econometric models break down when policy regimes change, emphasizing the importance of capturing invariant ``deep parameters.''.
    \item \textbf{Linear Assumptions}: \citet{tong1990non} and \citet{hamilton1989new} show that real markets exhibit nonlinear behaviors, threshold effects, and regime changes that linear models cannot adequately capture. \citet{terasvirta2005linear} demonstrates that nonlinear models often significantly outperform linear alternatives in financial time series forecasting.
    \item \textbf{Static Parameter Structure}: \citet{stock1996evidence} and \citet{pesaran2007selection} show that traditional models assume constant relationships between variables, ignoring temporal evolution and structural breaks.
\end{enumerate}

\subsection{Dynamical Systems Approach to Market Modeling}

Dynamical systems theory offers an alternative framework that explicitly models how variables interact and evolve according to mechanistic rules. This approach has enabled quantitative modeling across diverse fields:

\begin{itemize}
    \item \textbf{Population Dynamics}: As \citet{murray2002mathematical} shows, Lotka-Volterra equations model predator-prey relationships and competition between species
    \item \textbf{Epidemiology}: \citet{anderson1991infectious} demonstrate how SIR (Susceptible-Infected-Recovered) models track disease transmission through populations
    \item \textbf{Chemical Kinetics}: \citet{epstein1998introduction} describe how reaction-diffusion equations model concentration changes in chemical systems
    \item \textbf{Fluid Mechanics}: As \citet{tritton1988physical} explains, Navier-Stokes equations govern fluid flow and turbulence
\end{itemize}

Researchers have adapted these principles to economic and financial systems (see Appendix B). The key advantage of this approach is its explicit representation of how market variables change as functions of current conditions and external shocks. This mechanistic perspective enables:

\begin{enumerate}
    \item \textbf{Causal Understanding}: As \citet{strogatz2014nonlinear} demonstrates, ODEs can incorporate known economic relationships and behavioral assumptions, providing interpretable parameters that represent real-world processes.
    \item \textbf{Intervention Analysis}: By modifying parameters or adding forcing terms, \citet{blanchard1980solution} show how researchers can simulate the effects of policy changes, market interventions, or external shocks.
    \item \textbf{Nonlinear Dynamics}: As \citet{medio2001nonlinear} explains, ODE models naturally accommodate nonlinear relationships, feedback loops, and threshold effects commonly observed in markets.
    \item \textbf{Stability Analysis}: \citet{wiggins2003introduction} describes how mathematical techniques from dynamical systems theory can identify equilibrium points, assess stability, and predict long-term market behavior.
\end{enumerate}

\subsubsection{Ordinary Differential Equations in Market Context}

In the proposed framework, market variables such as prices, demand, supply, and market share are treated as state variables that evolve continuously over time. The rates of change of these variables are expressed as functions of the current market state, external economic factors, and model parameters that capture fundamental market mechanisms.

Several researchers have laid the groundwork for this approach:

\textbf{Price Dynamics}: \citet{zeeman1974unstable} developed catastrophe theory models of stock market crashes, while \citet{day1990bulls} analyzed nonlinear price adjustment mechanisms in commodity markets.

\textbf{Supply and Demand Dynamics}: \citet{saari1985iterative} studied the stability of market equilibria using differential equation methods, while \citet{scarf1960some} analyzed price adjustment processes in general equilibrium.

\textbf{Market Share Competition}: \citet{hanssens2001market} applied differential equations to model competitive dynamics in marketing, while \citet{chintagunta1992empirical} developed continuous-time models of brand choice.


\subsection{Parameter Estimation and Model Validation}

The parameter estimation process involves constrained optimization techniques that fit the ODE model to observed market data while ensuring that parameter values remain within economically meaningful ranges. This approach has been successfully applied in various contexts:

\textbf{Econometric Methods}: \citet{bergstrom1990continuous} developed continuous-time econometric methods for estimating differential equation models from discrete data, while \citet{phillips1991optimal} analyzed the statistical properties of such estimators.

\textbf{Optimization Techniques}: Recent advances in computational methods have made it feasible to estimate complex ODE models. \citet{ramsay2007parameter} developed parameter cascading methods, while \citet{brunel2008parameter} proposed likelihood-based approaches for stochastic differential equations.

\textbf{Model Selection and Validation}: \citet{burnham2002model} provide frameworks for model selection in complex systems, while cross-validation techniques adapted for time series by \citet{bergmeir2012use} can assess out-of-sample performance.

This research aims to demonstrate that ODE-based market models can provide superior performance in both forecasting accuracy and decision support capabilities compared to traditional AR approaches, while offering deeper insights into the fundamental processes governing market dynamics. The approach builds on the rich literature in complexity economics while leveraging modern computational methods to make such models practically implementable for real-world market analysis.

\section{Market Model Set-up}

\subsection{Definitions}
Let $M$ denote the market of $n \in \mathbb{N}$ products $d$. At any given time $t$, product $d \in M$ belongs to exactly one market segment: $i$ and one use case: $u$.\\
\\
\noindent At any time $t$, a customer segment $D_{i}(t)$ is the count of all products in the segment. A segment can be defined as any grouping of products---by brand, year of release, feature set, etc.\\
\\
\noindent At any time $t$, a product has one of two deterministic states: \textbf{active} or \textbf{obsolete}, which implies that $ND$ and $OD$ are mutually exclusive groups at time $t$.\\
\\
\noindent At any time $t$, for a given customer segment $i \in \{1,n\}$ where $\{1,n\} \subseteq \mathbb{N}$, the size of a market segment is given by its sources minus sinks:
\begin{equation}
\frac {d(D_i(t))}{dt} = \frac{d(ND_i(t))}{dt} - \frac{d(OD_i(t))}{dt}
\end{equation}

\indent Where, \\
\indent $D$: Total products\\
\indent $ND$: New products\\
\indent $OD$: Obsolete products\\

\begin{equation}
\frac {d(ND_i(t))}{dt} = \frac{d(BND_i(t))}{dt} + \frac{d(RD_i(t))}{dt}
\end{equation}

\indent Where, \\
\indent $ND$: New products\\
\indent $BND$: De novo products\\
\indent $RD$: Refreshed products\\

\noindent Putting these equations together, we get: 

\begin{equation}
\frac {d(D_i(t))}{dt} =\left(\frac{d(BND_i(t))}{dt} + \frac{d(RD_i(t))}{dt}\right) - 
\frac{d(OD_i(t))}{dt}
\end{equation}

$ND$ represents all demand entering a customer segment $i$, which constitutes demand due to the replacement of other segments $OD_{j}$ or $RD_i$ and fresh demand as a result of new entrants in the market $BND_i$. The proposed product market model exhibits structural similarities to the classical Lotka-Volterra (LV) predator-prey system, though with important conceptual differences that reflect the unique nature of technology markets (see Appendix A).

\section{Market Dynamics Components}
We now specify three components that govern segment evolution: obsolescence, refresh, and new entrant dynamics. We use the competitiveness scores of segments and helper functions described in sections below to allocate demand, modify competition scores, incorporate customer psychology, etc.

\subsection{Obsolescence Dynamics}

A segment becomes obsolete at a rate that decreases with its competitiveness and with customer attachment (stickiness):

\begin{equation}
\frac {d(OD_i(t))}{dt} = -\lambda \cdot D_{i}(t)
\text{, where }
\lambda = -k \cdot (1 - a_{i,m}(t)) \cdot (1 - s)
\end{equation}

\indent $s \in [0,1]$: Stickiness or attachment between customers and their customer segments, i.e. customers' resistance to change\\
\indent $k \in [0,50]$: Default rate of decay when $s$ = 0 and $a_{i,m}(t)$ = 0. 
\\

Thus, more competitive segments experience slower obsolescence. Conversely, less competitive segments deplete more rapidly. The stickiness parameter $s$ captures customer attachment that operates independently of measured competitiveness.

\subsection{Refresh Dynamics}

\noindent When products become obsolete, we assume customers replace them with new products. As products in segment $i$ become obsolete, demand flows to other segments according to their relative competitiveness. We allocate this demand using the mechanisms described in previous sections.

\begin{equation}
\frac{d}{dt}RD_i(t) = h(a_{i,m}^{mod}(t))\sum_{j = 1}^{n} \frac{d(OD_{j}(t))}{dt}
\end{equation}

Existing customers face the same choices as new customers when replacing obsolete products. The simplest assumption allocates all refresh demand to the most competitive segment $i=i_{max}$, corresponding to a winner-take-all scenario ($wta = 1$). This assumption gives,

\begin{equation}
\frac{d}{dt}RD_i(t) = \delta_{i,i_{max}} \cdot \sum_{j = 1}^{n} \frac{d(OD_{j}(t))}{dt} \text{, where $\delta_{m,n}$ is the Kronecker delta}
\end{equation}

\noindent Equivalently, this can be written using matrix notation as:
\begin{equation*}
\begin{bmatrix}
\frac{d(RD_1(t))}{dt} \\\\ \frac{d(RD_{2}(t))}{dt} \\ . \\ . \\ \frac{d(RD_{n-1}(t))}{dt} \\\\ \frac{d(RD_{n}(t))}{dt}
\end{bmatrix} = 
\sum_{j = 1}^{n} X_{i,j}
\end{equation*}
Where,
\begin{equation*}
X_{i,j} =
\begin{bmatrix}
\delta_{1,i_{max}} \\\\  
\delta_{2,i_{max}} \\ . \\ . \\
\delta_{n-1,i_{max}} \\\\ 
\delta_{n,i_{max}} 
\end{bmatrix}
\cdot 
\begin{bmatrix}
\frac{d(OD_1(t))}{dt} \\\\ \frac{d(OD_{2}(t))}{dt} \\ . \\ . \\ \frac{d(OD_{n-1}(t))}{dt} \\\\ \frac{d(OD_{n}(t))}{dt}
\end{bmatrix}^T
\end{equation*}

\noindent This winner-take-all (WTA) assumption can be relaxed, as shown previously, by allowing the allocation of demand to vary between a WTA scenario and a scenario in which demand is allocated as a function of competitiveness. In this case, the Kronecker delta vector (that necessarily sums to 1) is replaced with another vector based on competitiveness scores (that also necessarily sums to 1). Note that the competitiveness scores can be modified according to psychology and winner-take-all factors.

\begin{equation*}
\begin{bmatrix}
\frac{d(RD_1(t))}{dt} \\\\ \frac{d(RD_{2}(t))}{dt} \\ . \\ . \\ \frac{d(RD_{n-1}(t))}{dt} \\\\ \frac{d(RD_{n}(t))}{dt}
\end{bmatrix} = 
\begin{bmatrix}
h(a_{1,m}^{mod}(t)) \\\\  
h(a_{2,m}^{mod}(t)) \\ . \\ . \\
h(a_{n-1,m}^{mod}(t)) \\\\ 
h(a_{n,m}^{mod}(t)) 
\end{bmatrix}
\cdot 
\begin{bmatrix}
\frac{d(OD_1(t))}{dt} \\\\ \frac{d(OD_{2}(t))}{dt} \\ . \\ . \\ \frac{d(OD_{n-1}(t))}{dt} \\\\ \frac{d(OD_{n}(t))}{dt}
\end{bmatrix}^T
\end{equation*}
Where $h$ is a normalizing operation, or a demand allocation function, as described above.\\

\noindent To this point, we have not allowed for the possibility that a customer's refresh purchase decision might depend on the product that is becoming obsolete. There are many reasons to believe that such dynamics could be important. For example, a familiarity bias would necessitate a refresh demand equation that increases the allocation of demand to the product currently held, regardless of competitiveness scores. These dynamics require an $n\times n$ matrix rather than an $n$-length vector.

\begin{equation*}
\begin{bmatrix}
\frac{d(RD_1(t))}{dt} \\\\ \frac{d(RD_{2}(t))}{dt} \\ . \\ . \\ \frac{d(RD_{n-1}(t))}{dt} \\\\ \frac{d(RD_{n}(t))}{dt}
\end{bmatrix} = 
\begin{bmatrix} 
h(a_{1,1}^{mod}(t)) & h(a_{1,2}^{mod}(t)) & . & h(a_{1,n}^{mod}(t)) \\\\ 
h(a_{2,1}^{mod}(t)) & h(a_{2,2}^{mod}(t)) & . & h(a_{2,n}^{mod}(t)) \\
. & . & . & .\\  
. & . & . & .\\ 
h(a_{n-1,1}^{mod}(t)) & h(a_{n-1,2}^{mod}(t)) & . & h(a_{n-1,n}^{mod}(t)) \\\\
h(a_{n,1}^{mod}(t)) & h(a_{n,2}^{mod}(t)) & . & h(a_{n,n}^{mod}(t)) \\ 
\end{bmatrix} 
\cdot 
\begin{bmatrix}
\frac{d(OD_{1}(t))}{dt} \\\\  \frac{d(OD_{2}(t))}{dt} \\ . \\ . \\ \frac{d(OD_{n-1}(t))}{dt} \\\\ \frac{d(OD_{n}(t))}{dt}
\end{bmatrix} 
\end{equation*}

\noindent Here, as before, the $h$ operation must be defined such that each column of the allocation matrix sums to 1, i.e.
\begin{equation}
\sum_{j=1}^{n} h(a_{i,j}(t)) = 1
\end{equation}

\subsection{New Market Entrants}

A new consumer in the market chooses to purchase a product based on how competitive it is in the market. Similar to product refresh dynamics, new customers may either take a winner-take-all (WTA) approach or a more relaxed one.\\
\\
If we assume all new demand is allocated to the most competitive product $i=i_{max}$, the rate at which segment $D_{i}$ grows due to these de novo purchases is: 

\begin{equation}
\frac {d(BND_i(t))}{dt} = \delta_{i,i_{max}} \cdot \frac{d(NC(t))}{dt}
\end{equation}

\noindent Where, $\delta_{m,n}$ is the Kronecker delta and $\frac{d(NC(t))}{dt}$ is rate at which new customers are entering the market agnostic to specific customer segments at a time $t$.\\
\\
\noindent Equivalently, this can be written using matrix notation as:

\begin{equation*}
\begin{bmatrix}
\frac{d(BND_1(t))}{dt} \\\\ \frac{d(BND_{2}(t))}{dt} \\ . \\ . \\ \frac{d(BND_{n-1}(t))}{dt} \\\\ \frac{d(BND_{n}(t))}{dt}
\end{bmatrix} = 
\begin{bmatrix}
\delta_{1,i_{max}} \\\\  
\delta_{2,i_{max}} \\ . \\ . \\
\delta_{n-1,i_{max}} \\\\ 
\delta_{n,i_{max}} 
\end{bmatrix}
\cdot
\frac{d(NC(t))}{dt}
\end{equation*}

\noindent Identical to customer choices in the case of product refresh, the WTA assumption can be relaxed, by allowing the allocation of demand to vary between a WTA scenario, and a scenario in which demand is allocated as a function of competitiveness. In this case, the Kronecker delta vector (that necessarily sums to 1) is replaced with another vector based on competitiveness scores (that also necessarily sums to 1).

\begin{equation}
\frac {d(BND_i(t))}{dt} = h(a_{i,m}^{mod}(t)) \cdot \frac{d(NC(t))}{dt}
\end{equation}

\noindent Equivalently, this can be written using matrix notation as:
\begin{equation*}
\begin{bmatrix}
\frac{d(BND_1(t))}{dt} \\\\ \frac{d(BND_{2}(t))}{dt} \\ . \\ . \\ \frac{d(BND_{n-1}(t))}{dt} \\\\ \frac{d(BND_{n}(t))}{dt}
\end{bmatrix} = 
\begin{bmatrix}
h(a_{1,m}^{mod}(t)) \\\\  
h(a_{2,m}^{mod}(t)) \\ . \\ . \\
h(a_{n-1,m}^{mod}(t)) \\\\ 
h(a_{n,m}^{mod}(t)) 
\end{bmatrix}
\cdot
\frac{d(NC(t))}{dt}
\end{equation*}

\section{Competitiveness Scores}
The performance and competitiveness of a market segment is arguably the most pertinent factor that governs its growth or decline. We quantify various traits and features of market segments and define two measures: a pairwise and market-level score that provides some insight into the standing of a market segment in the broader landscape.\\
\\
Given any two customer segments $D_i, D_j \in M$ with $n$ segments, $a_{i,j}(t)$ is the \textbf{pairwise competitiveness score} of $D_i$ w.r.t $D_j$ at time $t$ and is a function attribute of the two market segments. 
\\
Furthermore, $a_{i,m}(t)$ is the \textbf{relative competitiveness score} of $D_i$ w.r.t to the rest of the customer segments in the market $M$, and is calculated by combining all of the pairwise competitiveness scores of $D_{i}$ at time $t$. 

\subsection{Scoring Methodology}
Customers score various product attributes, for example the price and quality of a segment, and indicate how important these attributes are to them.\\
\\
We state that there are $k$ number of distinguishable product attributes $z$, $count(z) = k$, that buyers evaluate on a predefined numerical scale: [$s_{min}, s_{max}] > 0$ independent of $t$:
\begin{equation}
z_{i}(t), z^{imp}_{i}(t) \in [s_{min},s_{max}]
\end{equation}
Here, $z_{i}(t)$ is the performance score of an attribute $z_{i}$ for segment $i$, and $z^{imp}_{i}(t)$ is the importance of the attribute to a customer. Notably, the performance score of an attribute is a raw score that reflects how a consumer values a feature, and is distinct from a customer segment's competitiveness score, which is defined later in this section. \\
\\
The "applied" importance score or \textit{weight} of an attribute $z$ is the relative importance of said attribute to all other attributes that a customer is evaluating:
\begin{equation}
w^{z, imp}_{i}(t) = \frac{z^{imp}_{i}(t)}{\sum_{z=1}^{k} z^{imp}_{i}(t)} \text{, where } w^{z, imp}_{i}(t) \in [0,1]
\end{equation}

The relative competitiveness score of a customer segment $i$ in the market of products at time $t$ is:
\begin{equation}
a_{i,m}(t) = \sum_{z = 1}^{k} z_{i}(t) \cdot w^{z, imp}_{i}(t) \text{, where }
a_{i,m}(t) \in [s_{min},s_{max}]
\end{equation}

\textbf{Example}: If price and quality are two attributes a customer values in their decision-making, $k = 2$, $z_{1} = p$ and $z_{2} = q$, and $[s_{min},s_{max}] = [0,10]$, the competitiveness score of segment $i$ compared to the remaining segments in the market is:
\begin{equation}
a_{i,m}(t) = [p_{i}(t), q_{i}(t)] \cdot \left[ \frac{p^{imp}_{i}(t)}{p^{imp}_{i}(t) + q^{imp}_{i}(t)}, \frac{q^{imp}_{i}(t)}{p^{imp}_{i}(t) + q^{imp}_{i}(t)} \right]\text{, where }
a_{i,m}(t) \in [1,10]
\end{equation}

Next, the pairwise competitiveness score between segment $i$ and $j$ is:
\begin{equation}
a_{i,j}(t) =
\sum_{z = 1}^{k}
\frac{z_{i}(t) - z_{j}(t)}{s_{max} - s_{min}} \cdot w^{z, imp}_{i}(t)\text{, where }
a_{i,j}(t) \in [-1,1]
\end{equation}
Here, $\frac{z_{i}(t) - z_{j}(t)}{s_{max} - s_{min}}$ is the relative performance score of attribute $z$ in segment $i$ and $j$. Conceptually, this is asking the following question: ``How much more competitive is $i$ than $j$ relative to how much more competitive $i$ \emph{could be} than $j$.''\\
\\
In simple terms, $a_{i,j}(t)$ is the degree to which customer segment $i$ is better (or worse) than customer segment $j$. A negative score means that segment $i$ is worse than $j$, a score of 0 means that neither segment is better than the other, and a positive score means that $i$ is better than segment $j$.\\

Consequently,
\begin{equation}
a_{i,j}(t) = -a_{j,i}(t)
\end{equation}

\section{Market Building Blocks}
Now that we have established the core parameters that drive market dynamics, i.e. the flow of demand through the network of market segments, we define various modular functions that can be applied in a plug-and-play fashion to model the growth and decline of segment sizes. 

\subsection{Behavioral Modification}
\subsubsection{Winner-Take-All Effects}
The competitiveness score of a segment is inflated or deflated based on how winner-take-all the market is, i.e. how logical or rational customers are (at large) by choosing to buy products that are in fact the best for them.\\
\\
We already defined the pairwise competitiveness score, $a_{i,j}(t)\in [-1,1]$ and the market level competitiveness score, $a_{i,m}(t) \in [s_{min}, s_{max}]$, where $s_{min} \geq 0$\\
\\
The winner-take-all factor is defined at the market level and is not dependent on the individual competitiveness scores:
\begin{equation}
wta \in [0,1]
\end{equation}
Say the most competitive product in the market is $i=i_{max}$ with a score of $a_{m (j)}^{max}(t)$.\\
\\
The competitiveness scores are inflated or deflated as follows using the $wta$ factor:
\begin{equation}
a_{i,m (j)}^{mod}(t) =
    \begin{cases}
    (\delta_{i,i_{max}} \cdot a_{i,m (j)}(t) \cdot wta) + (a_{i,m (j)}(t) \cdot (1 - wta)), & \text{if $a_{i,m (j)}(t) > 0$}.\\
    a_{i,m (j)}(t) - |a_{i,m (j)}(t) \cdot wta|, & \text{if $a_{i,m (j)}(t) < 0$}.\\
    - 1 \cdot wta, & \text{if $a_{i,m (j)}(t) = 0$}.
    \end{cases}
\end{equation}
\noindent Where $\delta_{m,n}$ is the Kronecker delta function.
\\
In essence, we discard $wta$\% of the competitiveness score in all cases except the first case when $i = i_{max}$, which implies that $a_{i,m (j)}^{mod}(t) = a_{i,m (j)}(t) = a_{m (j)}^{max}(t)$. If $wta < 1$ in the first case, then for  $i \neq i_{max}$, the logic of discarding $wta$\% of the competitiveness score stands since $\delta_{i,i_{max}} = 0$.\\
\\
If $wta = 1$ in the first case, $a_{i,m (j)}^{mod}(t)$ = $\delta_{i,i_{max}} \cdot a_{i,m (j)}(t)$. This implies that the competitiveness score of the best product in the market remains unchanged while the scores of the remaining products are modified to 0.

\subsubsection{Customer Psychology}
The true competitiveness of a segment may vary slightly from how a customer applies this competitiveness knowledge in their decisions. For example, a product may be highly competitive based on the attributes that a customer values, but if it is a new unknown brand (e.g., ASUS) or is manufactured by a non-European country (e.g., Egypt), they may value the product less than it scored. We assert that this behavior can be quantified as a difference between a \textit{true competitiveness score} and an \textit{applied competitiveness score}.\\
\\
We assert that customer resistance to rationality decreases as the product competitiveness score increases in the market. Hence, we model this resistance as an exponential curve whose rate of growth (or decline) is proportional to true market competitiveness. We assert that a customer's perception of competitiveness cannot become skewed in the opposite direction, i.e. consider a product that is worse than what they currently own as better.\\
\\
We already defined the pairwise competitiveness score, $a_{i,j}(t)\in [-1,1]$ and the market level competitiveness score, $a_{i,m}(t) \in [s_{min}, s_{max}]$, where $s_{min} \geq 0$ and $s_{max} \geq s_{min}$\\
\\
We define the following resistance metrics to quantify the deviation from the true competitiveness score. 
\begin{equation}
    r(a_{i,j}(t), \gamma, c) = 
    \begin{cases}
        1, & \text{if $a_{i,j}(t) \leq 0$}.\\
        e^{\gamma \cdot a_{i,j}(t)} + c, & \text{if $a_{i,j}(t) > 0$}
    \end{cases}
\end{equation}
\begin{equation}
r(a_{i,m}(t), \gamma, c) = 
    \begin{cases}
        1, & \text{if $a_{i,m}(t) = 0$}.\\
        e^{\gamma \cdot a_{i,m}(t)} + c, & \text{if $a_{i,m}(t) > 0$}.
    \end{cases}
\end{equation}
\begin{equation}
\text{Where, } \gamma < 0, c \geq 0 
\text{ and }
r \in [0,1]
\end{equation}
$c$ represents the default degree of resistance a customer has towards buying a segment, even the most competitive one in the market.\\
\\
We discard a portion of the competitiveness score in the same way as we deal with the winner-take-all characteristic of the market. The competitiveness scores are modified as follows:
\begin{equation}
a_{i,m (j)}^{mod}(t) =
    \begin{cases}
    a_{i,m (j)}(t) \cdot (1 - r), & \text{if $a_{i,m (j)}(t)  > 0$}.\\
    a_{i,m (j)}(t) - |a_{i,m (j)}(t) \cdot r|, & \text{if $a_{i,m (j)}(t) < 0$}.\\
    0, & \text{if $a_{i,m (j)}(t) = 0$}.
    \end{cases}
\end{equation}

\section{Demand Allocation Mechanisms}
In order to operationalize our competitiveness score in market share allocation, we define a function $h(a_{i,m (j)}(t))$ distributes demand across all segments in the market.

The distinction in the approaches described below primarily hinge on same-product refreshes---bias towards it, bias against it, or somewhere in between. These biases and preferences reflect different assumptions about the market. 

\subsection{Ratio-Based Allocation}

Demand is allocated for positive competitiveness scores only, i.e. products that are indeed better than the existing segment $i$. Notably, this allocation method does not return any market share to the existing segment $i$ and instead increases demand in segments that are positive improvements from the focus segment.\\
$\forall i \in n$ and $\forall j \in n,$ if $a_{i,m (j)}(t) \geq 0$:

\begin{equation}
h(a_{i,m (j)}(t)) = \frac {a_{i,m (j)}(t)}{\sum_{i=1}^{i=n} a_{i,m (j)}(t)}, \text{where } {a_{i,m (j)}(t)} \geq 0
\end{equation}

If $a_{i,m (j)}(t) < 0$, no demand is allocated to the less competitive segment or product.
\begin{equation}
h(a_{i,m (j)}(t)) = 0
\end{equation}
Consequently,
\begin{equation}
\sum_{i=1}^{i=n} h(a_{i,m (j)}(t)) = 1
\end{equation}

\subsection{Softmax-Based Allocation}
Alternatively, we can distribute market demand using a standard softmax function. We apply softmax on the vector of segment competitiveness scores and index the resultant vector at position $i$ for the final allocation towards segment $i$ in the market.

\begin{equation}
h(a_{i,m (j)}(t)) = softmax( 
 \left\langle a_{i,m (j)}(t), a_{i+1,m (j)}(t) ... a_{n,m (j)}(t) \right\rangle) [i]
\end{equation}
\\
By nature of the softmax algorithm,
\begin{equation}
\sum_{i=1}^{i=n} h(a_{i,m (j)}(t)) = 1
\end{equation}
Notably, in contrast to \textbf{Option 1}, negative competitiveness scores, i.e. products segments that are worse than the existing one, still win a non-trivial portion of the demand. 

\subsection{Redistribution Allocation}

It is evident from both \textbf{Option 1} and \textbf{Option 2} that market demand can be allocated towards segments that are less competitive than the existing one in a multitude of ways. In \textbf{Option 1}, no demand is allocated to a poorer performing segment, which is modified in \textbf{Option 2}. \\
\\
Here, similar to \textbf{Option 1} no demand is allocated towards poorer performing segments. If $a_{i,m (j)}(t) < 0$, 
\begin{equation}
h(a_{i,m (j)}(t)) = 0
\end{equation}\\
However, in contrast to \textbf{Option 1}, instead of allocating additional demand across the stronger segments, we allocate remaining demand towards the existing segment itself. 

\begin{equation}
h(a_{i,m (j)}(t)) = a_{i,m (j)}(t) \text{ if $a_{i,m (j)}(t)  > 0$}.
\end{equation}
\begin{equation}
h(a_{i,i}(t)) = 1 - \sum a_{i,m (j)}(t) \text{ where } {a_{i,m (j)}(t) > 0}
\end{equation}

In a winner-take-all scenario, the WTA factor defined in the previous section is set to 1. For all allocation methodologies---ratio, softmax and redistribution, we arrive at the same outcome:

\begin{equation}
h(a_{i,j(m)}(t)) = \delta_{i,i_{max}} \cdot a_{i,m (j)}(t) = 
\begin{cases}
    1 \text{ where } {i = i_{max}}\\
    0 \text{ where } {i \neq i_{max}}
\end{cases}
\end{equation}
In essence, all of the demand is allocated towards the most competitive segment in the market $i_{max}$.

\section{Sample Implementation}

Consider a market of $n = 3$ customer segments $D_{i}$ and $t = 6$ discrete time stamps. Each segment $i$ earns the following raw performance scores: 

\begin{table}[H]
\centering
\caption{Quality and Price Scores for Segments Over Time}
\begin{tabular}{c|ccc|ccc}
\toprule
& \multicolumn{3}{c|}{Quality (scores $\in$ [1, 10])} & \multicolumn{3}{c}{Price (scores $\in$ [1, 10])} \\
Time & $D_{i = 1}$ & $D_{i = 2}$ & $D_{i = 3}$ & $D_{i = 1}$ & $D_{i = 2}$ & $D_{i = 3}$ \\
\midrule
$t_1$ & 4 & 8 & 2 & 8 & 4 & 8 \\
$t_2$ & 5 & 8 & 2 & 8 & 5 & 8 \\
$t_3$ & 6 & 7 & 2 & 7 & 6 & 8 \\
$t_4$ & 7 & 6 & 2 & 6 & 7 & 8 \\
$t_5$ & 8 & 5 & 3 & 5 & 8 & 8 \\
$t_6$ & 8 & 4 & 3 & 4 & 8 & 8 \\
\bottomrule
\end{tabular}
\end{table}
Note that some segments earn better performance scores for certain attributes over time while others receive lower scores. This reflects evolutions in customer preferences of segment attributes. The competitiveness scores of each segment is a combination of the directionality of attribute-specific performance scores. Over time, the segment performance and competitiveness scores follow the trajectory below:

\begin{figure}
    \centering
    \begin{subfigure}[b]{0.45\textwidth}
        \centering
        \includegraphics[width=\textwidth]{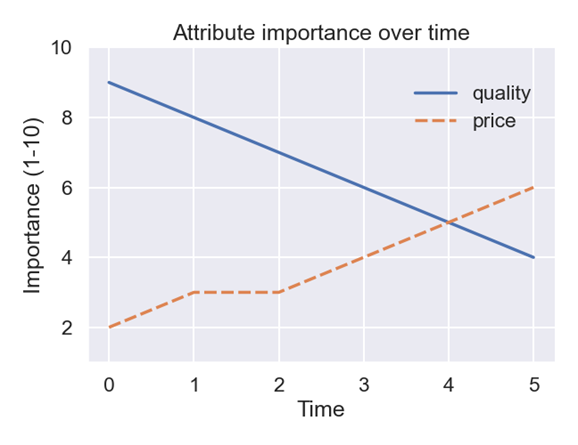}
        \caption{Importance Scores Over Time}
        \label{fig:imp_score}
    \end{subfigure}
    \begin{subfigure}[b]{0.45\textwidth}
        \centering
        \includegraphics[width=\textwidth]{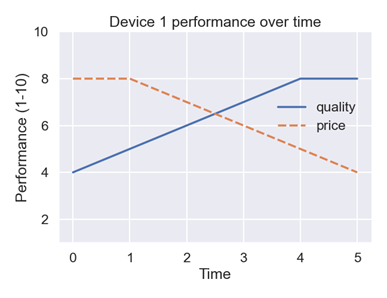}
        \caption{Segment 1 Performance Over Time}
        \label{fig:d1_perf}
    \end{subfigure}
    
    \begin{subfigure}[b]{0.45\textwidth}
        \centering
        \includegraphics[width=\textwidth]{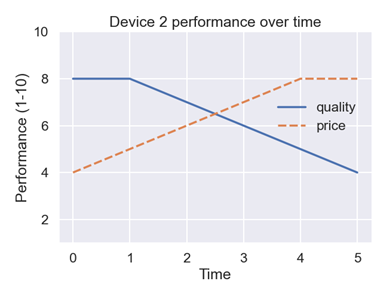}
        \caption{Segment 2 Performance Over Time}
        \label{fig:d2_perf}
    \end{subfigure}
    \begin{subfigure}[b]{0.45\textwidth}
        \centering
        \includegraphics[width=\textwidth]{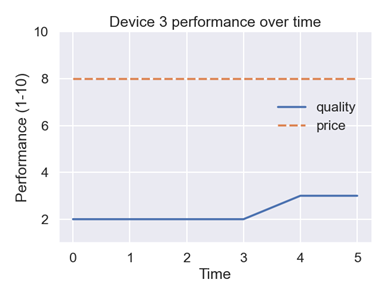}
        \caption{Segment 3 Performance Over Time}
        \label{fig:d3_perf}
    \end{subfigure}
    
    \caption{Segment Performance Analysis Over Time}
    \label{fig:product_analysis}
\end{figure}

We see that over time $t=0 \text{ to } 6$, $D_3$ becomes more competitive in the market while $D_2$ becomes less competitive. Ultimately, at $t = 6$, $D_2$ is the most competitive customer segment in the market. Next, demand is allocated across the three market segments at each time stamp using one of three methodologies---ratio, softmax or redistribution. We also incorporate a winner-take-all factor of 0.3. 

\begin{figure}[H]
    \centering
    \includegraphics[width=0.5\linewidth]{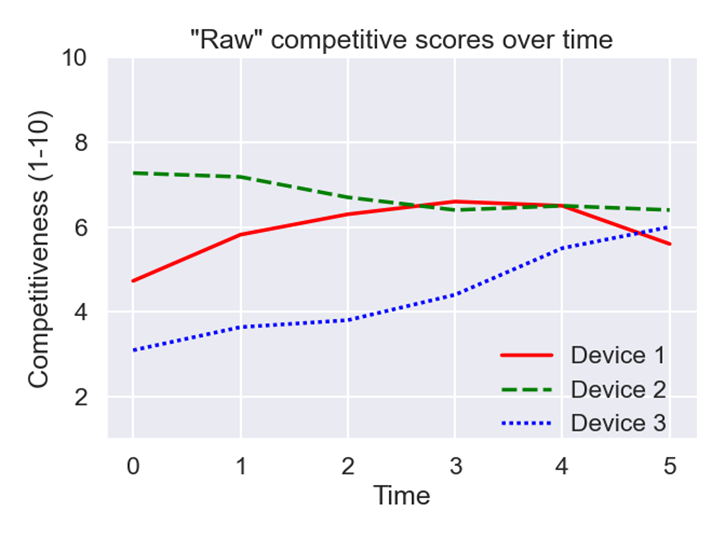}
    \caption{Segment Competitiveness Over Time}
    \label{fig:product_competitiveness}
\end{figure}

\begin{figure}[H]
    \centering
    \begin{subfigure}[b]{0.45\textwidth}
        \centering
        \includegraphics[width=\textwidth]{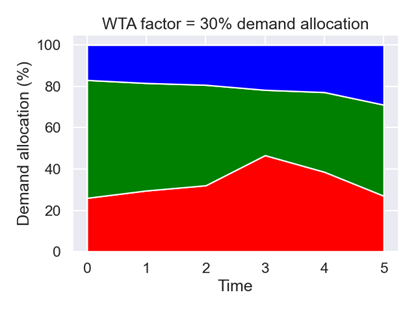}
        \caption{Moderate winner-take-all ($wta = 0.3$)}
        \label{fig:wta_30}
    \end{subfigure}
    \begin{subfigure}[b]{0.45\textwidth}
        \centering
        \includegraphics[width=\textwidth]{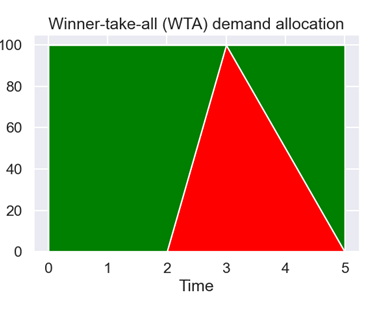}
        \caption{Complete winner-take-all ($wta = 1$)}
        \label{fig:wta_1}
    \end{subfigure}
    
    \caption{Redistribution-based Demand Allocation Under Different Market Scenarios}
    \label{fig:wta_comparison}
\end{figure}

We see that in the scenario with $wta=1$, the most competitive device in the market at any $t$ receives 100\% of the market demand. In the scenario with $wta = 0.3$, we interpolate between a scenario of $wta = 1$ and $wta = 0$, leaning more closely to the latter. Demand is therefore distributed more proportionally across the three market segments.

\section{Future Work}

Real-world use cases are likely to have historical market shares and performance estimates, but not customer importance scores or WTA factors. These parameters can be discovered and tuned with constrained optimization. Iteratively collecting and comparing market share data with model predictions can offer important insights. Deviations from model outputs are important learning opportunities and can provide a remedy for speculative forecasting.

\section{Conclusion}

This paper presents a novel framework for modeling market dynamics that explicitly incorporates product competitiveness and consumer behavior. The key contributions include: (1) a mechanistic modeling framework that moves beyond traditional AR approaches by explicitly modeling underlying market processes, (2) a comprehensive competitiveness scoring system that quantifies product attributes and consumer preferences, and (3) flexible demand allocation mechanisms accommodating different market behaviors from winner-take-all to distributed competitive outcomes.

Our framework addresses fundamental limitations of traditional econometric models by providing causal understanding, enabling scenario analysis, naturally accommodating nonlinear dynamics, and offering interpretable parameters representing real-world processes. The sample implementation demonstrates how competitiveness affects market share evolution, while behavioral factors such as customer psychology provide more nuanced understanding than purely statistical approaches.

By bridging concepts from mathematical biology, dynamical systems theory, and economics, this framework represents a significant step toward more mechanistic, interpretable, and actionable market models that can support strategic decision-making in complex competitive environments.

\section{Acknowledgment}

We would like to thank our colleagues and management at Intel for providing business use cases and problem statements that motivated our modeling approaches. We are also grateful to the Institute for Operations Research and Management Science (INFORMS) for the opportunity to share our work with the larger OR community. 

\appendix

\section{Relationship to Lotka-Volterra Dynamics}

The standard LV equations describe the dynamics of two interacting populations:

\begin{align}
\frac{dx}{dt} &= \alpha x - \beta xy \\
\frac{dy}{dt} &= \delta xy - \gamma y
\end{align}

Where $x$ represents prey population, $y$ represents predator population, and the parameters $\alpha, \beta, \gamma, \delta > 0$ represent birth rates, predation rates, and death rates.

In our product market model, equation (4) can be rewritten to highlight the analogous structure:

\begin{equation}
\frac{d(D_i(t))}{dt} = \underbrace{\frac{d(BND_i(t))}{dt}}_{\text{Birth term}} + \underbrace{\frac{d(RD_i(t))}{dt}}_{\text{Migration term}} - \underbrace{\frac{d(OD_i(t))}{dt}}_{\text{Death term}}
\end{equation}

The correspondence becomes clearer when we consider the competitive dynamics between market segments:

\subsection{Birth and Death Processes}
\begin{itemize}
    \item \textbf{Birth term} $\frac{d(BND_i(t))}{dt}$: Analogous to the prey birth rate $\alpha x$ in LV equations, representing the intrinsic growth of segment $i$ through new product introductions.
    
    \item \textbf{Death term} $\frac{d(OD_i(t))}{dt}$: Analogous to the predator death rate $\gamma y$, representing products becoming obsolete due to technological advancement or market forces.
\end{itemize}

\subsection{Competitive Interactions}
The \textbf{refresh term} $\frac{d(RD_i(t))}{dt}$ captures the most interesting parallel to LV dynamics. In technology markets, this term can represent:

\begin{itemize}
    \item \textbf{Competitive displacement}: When superior technology in segment $j$ causes users to abandon segment $i$, similar to predation in LV models
    \item \textbf{Market cannibalization}: When a company's new customer segment reduces demand for its existing segments
    \item \textbf{Technology migration}: Users moving between segments due to changing preferences or capabilities
\end{itemize}

\section{Applications of Dynamical Systems in Economics and Finance}

Dynamical systems theory has been successfully applied across various economic and financial domains, providing the foundation for our market modeling approach. 

\subsection{Macroeconomic Dynamics}

Business cycle modeling represents one of the earliest applications of biological dynamics to economics. As \citet{goodwin1967growth} showed, employment rates and wage shares interact cyclically in predator-prey fashion, with high employment leading to wage increases that eventually reduce profitability and employment. More recently, \citet{grasselli2017comment} have developed sophisticated stock-flow consistent models using differential equations that capture dynamic interactions between economic sectors while maintaining accounting identities.

\subsection{Behavioral Finance}

Another area ripe for dynamical systems applications is behavioral finance, where psychological factors and social interactions drive market dynamics. As \citet{lux1995herd} demonstrates, models of herd behavior explain how rational individual decisions can lead to collective market bubbles and crashes through purely endogenous mechanisms.\\
\\
These applications demonstrate key advantages of dynamical systems approaches: they model rates of change rather than levels, incorporate nonlinear relationships and feedback effects, provide mechanistic explanations for observed phenomena, and enable scenario analysis. Our market segmentation model builds on these foundations to address the specific challenges of modeling competitive dynamics in technology markets.

\bibliographystyle{plainnat}

\bibliography{references}

@book{arthur1997economy,
  title={The economy as an evolving complex system {II}},
  editor={Arthur, W. Brian and Durlauf, Steven N. and Lane, David A.},
  publisher={Addison-Wesley},
  address={Reading, MA},
  year={1997},
  url={https://www.santafe.edu/research/results/books/the-economy-as-an-evolving-complex-system-ii}
}

@article{farmer2009economy,
  title={The economy needs agent-based modelling},
  author={Farmer, J. Doyne and Foley, Duncan},
  journal={Nature},
  volume={460},
  number={7256},
  pages={685--686},
  year={2009},
  url={https://doi.org/10.1038/460685a}
}

@book{box1976time,
  title={Time series analysis: forecasting and control},
  author={Box, George E. P. and Jenkins, Gwilym M.},
  publisher={Holden-Day},
  address={San Francisco},
  year={1976}
}

@article{sims1980macroeconomics,
  title={Macroeconomics and reality},
  author={Sims, Christopher A.},
  journal={Econometrica},
  volume={48},
  number={1},
  pages={1--48},
  year={1980},
  url={https://doi.org/10.2307/1912017}
}

@article{granger1969investigating,
  title={Investigating causal relations by econometric models and cross-spectral methods},
  author={Granger, Clive W. J.},
  journal={Econometrica},
  volume={37},
  number={3},
  pages={424--438},
  year={1969},
  url={https://doi.org/10.2307/1912791}
}

@book{hamilton1994time,
  title={Time series analysis},
  author={Hamilton, James D.},
  publisher={Princeton University Press},
  address={Princeton, NJ},
  year={1994},
  url={https://press.princeton.edu/books/hardcover/9780691042893/time-series-analysis}
}

@book{pearl2009causality,
  title={Causality},
  author={Pearl, Judea},
  edition={2nd},
  publisher={Cambridge University Press},
  address={Cambridge},
  year={2009},
  url={https://doi.org/10.1017/CBO9780511803161}
}

@book{strogatz2014nonlinear,
  title={Nonlinear dynamics and chaos: with applications to physics, biology, chemistry, and engineering},
  author={Strogatz, Steven H.},
  edition={2nd},
  publisher={Westview Press},
  address={Boulder, CO},
  year={2014},
  url={https://www.westviewpress.com/books/nonlinear-dynamics-and-chaos}
}

@book{medio2001nonlinear,
  title={Nonlinear dynamics: a primer},
  author={Medio, Alfredo and Lines, Marji},
  publisher={Cambridge University Press},
  address={Cambridge},
  year={2001},
  url={https://doi.org/10.1017/CBO9780511754050}
}

@book{lutkepohl2005new,
  title={New introduction to multiple time series analysis},
  author={L{\"u}tkepohl, Helmut},
  publisher={Springer},
  address={Berlin},
  year={2005},
  url={https://doi.org/10.1007/978-3-540-27752-1}
}

@book{hendry1995dynamic,
  title={Dynamic econometrics},
  author={Hendry, David F.},
  publisher={Oxford University Press},
  address={Oxford},
  year={1995},
  url={https://global.oup.com/academic/product/dynamic-econometrics-9780198283164}
}

@article{lucas1976econometric,
  title={Econometric policy evaluation: {A} critique},
  author={Lucas, Robert E. Jr.},
  journal={Carnegie-Rochester Conference Series on Public Policy},
  volume={1},
  pages={19--46},
  year={1976},
  url={https://doi.org/10.1016/S0167-2231(76)80003-6}
}

@article{terasvirta2005linear,
  title={Linear models, smooth transition autoregressions, and neural networks for forecasting macroeconomic time series: {A} re-examination},
  author={Ter{\"a}svirta, Timo and Van Dijk, Dick and Medeiros, Marcelo C.},
  journal={International Journal of Forecasting},
  volume={21},
  number={4},
  pages={755--774},
  year={2005},
  url={https://doi.org/10.1016/j.ijforecast.2005.04.010}
}

@book{tong1990non,
  title={Non-linear time series: a dynamical system approach},
  author={Tong, Howell},
  publisher={Oxford University Press},
  address={Oxford},
  year={1990},
  url={https://global.oup.com/academic/product/non-linear-time-series-9780198522249}
}

@article{hamilton1989new,
  title={A new approach to the economic analysis of nonstationary time series and the business cycle},
  author={Hamilton, James D.},
  journal={Econometrica},
  volume={57},
  number={2},
  pages={357--384},
  year={1989},
  url={https://doi.org/10.2307/1912559}
}

@article{stock1996evidence,
  title={Evidence on structural instability in macroeconomic time series relations},
  author={Stock, James H. and Watson, Mark W.},
  journal={Journal of Business \& Economic Statistics},
  volume={14},
  number={1},
  pages={11--30},
  year={1996},
  url={https://doi.org/10.1080/07350015.1996.10524626}
}

@article{pesaran2007selection,
  title={Selection of estimation window in the presence of breaks},
  author={Pesaran, M. Hashem and Timmermann, Allan},
  journal={Journal of Econometrics},
  volume={137},
  number={1},
  pages={134--161},
  year={2007},
  url={https://doi.org/10.1016/j.jeconom.2006.03.010}
}

@book{murray2002mathematical,
  title={Mathematical biology: {I}. {A}n introduction},
  author={Murray, James D.},
  edition={3rd},
  publisher={Springer},
  address={New York},
  year={2002},
  url={https://doi.org/10.1007/b98868}
}

@book{anderson1991infectious,
  title={Infectious diseases of humans: dynamics and control},
  author={Anderson, Roy M. and May, Robert M.},
  publisher={Oxford University Press},
  address={Oxford},
  year={1991},
  url={https://global.oup.com/academic/product/infectious-diseases-of-humans-9780198545996}
}

@book{epstein1998introduction,
  title={An introduction to nonlinear chemical dynamics: oscillations, waves, patterns, and chaos},
  author={Epstein, Irving R. and Pojman, John A.},
  publisher={Oxford University Press},
  address={Oxford},
  year={1998},
  url={https://global.oup.com/academic/product/an-introduction-to-nonlinear-chemical-dynamics-9780195096705}
}

@book{tritton1988physical,
  title={Physical fluid dynamics},
  author={Tritton, D. J.},
  edition={2nd},
  publisher={Oxford University Press},
  address={Oxford},
  year={1988},
  url={https://global.oup.com/academic/product/physical-fluid-dynamics-9780198544937}
}

@article{blanchard1980solution,
  title={The solution of linear difference models under rational expectations},
  author={Blanchard, Olivier J. and Kahn, Charles M.},
  journal={Econometrica},
  volume={48},
  number={5},
  pages={1305--1311},
  year={1980},
  url={https://doi.org/10.2307/1912186}
}

@book{wiggins2003introduction,
  title={Introduction to applied nonlinear dynamical systems and chaos},
  author={Wiggins, Stephen},
  edition={2nd},
  publisher={Springer},
  address={New York},
  year={2003},
  url={https://doi.org/10.1007/b97481}
}

@article{zeeman1974unstable,
  title={On the unstable behaviour of stock exchanges},
  author={Zeeman, E. Christopher},
  journal={Journal of Mathematical Economics},
  volume={1},
  number={1},
  pages={39--49},
  year={1974},
  url={https://doi.org/10.1016/0304-4068(74)90034-2}
}

@article{day1990bulls,
  title={Bulls, bears and market sheep},
  author={Day, Richard H. and Huang, Weihong},
  journal={Journal of Economic Behavior \& Organization},
  volume={14},
  number={3},
  pages={299--329},
  year={1990},
  url={https://doi.org/10.1016/0167-2681(90)90061-H}
}

@article{saari1985iterative,
  title={Iterative price mechanisms},
  author={Saari, Donald G.},
  journal={Econometrica},
  volume={53},
  number={5},
  pages={1117--1131},
  year={1985},
  url={https://doi.org/10.2307/1911017}
}

@article{scarf1960some,
  title={Some examples of global instability of the competitive equilibrium},
  author={Scarf, Herbert},
  journal={International Economic Review},
  volume={1},
  number={3},
  pages={157--172},
  year={1960},
  url={https://doi.org/10.2307/2525561}
}

@book{hanssens2001market,
  title={Market response models: econometric and time series analysis},
  author={Hanssens, Dominique M. and Parsons, Leonard J. and Schultz, Randall L.},
  publisher={Springer},
  address={New York},
  year={2001},
  url={https://doi.org/10.1007/978-1-4615-1701-4}
}

@article{chintagunta1992empirical,
  title={An empirical investigation of advertising strategies in a dynamic duopoly},
  author={Chintagunta, Pradeep K. and Vilcassim, Naufel J.},
  journal={Management Science},
  volume={38},
  number={9},
  pages={1230--1244},
  year={1992},
  url={https://doi.org/10.1287/mnsc.38.9.1230}
}

@book{bergstrom1990continuous,
  title={Continuous time econometric modelling},
  author={Bergstrom, A. R.},
  publisher={Oxford University Press},
  address={Oxford},
  year={1990},
  url={https://global.oup.com/academic/product/continuous-time-econometric-modelling-9780198283409}
}

@article{phillips1991optimal,
  title={Optimal inference in cointegrated systems},
  author={Phillips, Peter C. B.},
  journal={Econometrica},
  volume={59},
  number={2},
  pages={283--306},
  year={1991},
  url={https://doi.org/10.2307/2938258}
}

@article{ramsay2007parameter,
  title={Parameter estimation for differential equations: a generalized smoothing approach},
  author={Ramsay, James O. and Hooker, Giles and Campbell, David and Cao, Jiguo},
  journal={Journal of the Royal Statistical Society: Series B (Statistical Methodology)},
  volume={69},
  number={5},
  pages={741--796},
  year={2007},
  url={https://doi.org/10.1111/j.1467-9868.2007.00610.x}
}

@article{brunel2008parameter,
  title={Parameter estimation of {ODE}'s via nonparametric estimators},
  author={Brunel, Nicolas J.-B.},
  journal={Electronic Journal of Statistics},
  volume={2},
  pages={1242--1267},
  year={2008},
  url={https://doi.org/10.1214/07-EJS132}
}

@book{burnham2002model,
  title={Model selection and multimodel inference: a practical information-theoretic approach},
  author={Burnham, Kenneth P. and Anderson, David R.},
  publisher={Springer},
  address={New York},
  year={2002},
  url={https://doi.org/10.1007/b97636}
}

@article{bergmeir2012use,
  title={On the use of cross-validation for time series predictor evaluation},
  author={Bergmeir, Christoph and Ben{\'\i}tez, Jos{\'e} M.},
  journal={Information Sciences},
  volume={191},
  pages={192--213},
  year={2012},
  url={https://doi.org/10.1016/j.ins.2011.12.028}
}

@incollection{goodwin1967growth,
  title={A growth cycle},
  author={Goodwin, Richard M.},
  booktitle={Socialism, Capitalism and Economic Growth},
  editor={Feinstein, C. H.},
  pages={54--58},
  publisher={Cambridge University Press},
  address={Cambridge},
  year={1967}
}

@article{grasselli2017comment,
  title={A comment on `{T}esting {G}oodwin: growth cycles in ten {OECD} countries'},
  author={Grasselli, Matheus and Maheshwari, Aditya},
  journal={Cambridge Journal of Economics},
  volume={41},
  number={6},
  pages={1761--1766},
  year={2017},
  url={https://doi.org/10.1093/cje/bex006}
}

@article{lux1995herd,
  title={Herd behaviour, bubbles and crashes},
  author={Lux, Thomas},
  journal={Economic Journal},
  volume={105},
  number={431},
  pages={881--896},
  year={1995},
  url={https://doi.org/10.2307/2235156}
}

\end{document}